\documentclass[a4paper]{article}
\usepackage[english]{babel}
\usepackage[utf8]{inputenc}
\usepackage{fullpage}
\usepackage{amsmath,amssymb,latexsym}
\usepackage[pdftex]{graphicx}
\usepackage{graphics}
\usepackage{psfrag}
\usepackage{color, epsfig}
\usepackage{mathrsfs}
\usepackage{theorem}  
\usepackage[colorlinks=true]{hyperref}

\def\R{{\mathbb R}}
\def\N{{\mathbb N}}

\def\S{{\mathcal S}}  

\def\ds{\displaystyle}

\def\div{\operatorname{div}\,} 
\def\divO{\operatorname{div}_{\partial\O}\,} 

\newcommand{\tr}{\operatorname{tr\,}}
\renewcommand{\leq}{\leqslant}
\renewcommand{\geq}{\geqslant}
\renewcommand{\O}{\mathcal{O}}
\renewcommand{\S}{\mathcal{S}}

\newtheorem{theorem}{Theorem}[section]

{\theorembodyfont{\rmfamily} \newtheorem{remark}[theorem]{Remark}}

\newcommand{\Finn}{\begin{flushright} \rule{2mm}{2mm}\end{flushright}}

\title{Computational Shape Derivatives in Heat Conduction: An Optimization Approach for Enhanced Thermal Performance
}

\author{
		M. Azaiez\thanks{Université de Bordeaux, I2M, (UMR CNRS 5295), 33607 Pessac, France, E-mail: {\tt azaiez@u-bordeaux.fr}.  Supported by the ANR project  NumOpTes ANR-22-CE46-0005 (2023--2026)} \and A. Doubova\thanks{Universidad de Sevilla, Dpto.\ EDAN e IMUS, Campus Reina Mercedes, 41012~Sevilla, Spain, E-mail: {\tt doubova@us.es}. Partially supported by PID2020-114976GB-I00, Plan estatal 2021, PY20 01125, PAIDI 2021.} \and S. Ervedoza\thanks{Institut de Mathématiques de Bordeaux, UMR 5251, Université de Bordeaux, CNRS, Bordeaux INP, F-33400 Talence, France,  
\texttt{sylvain.ervedoza@math.u-bordeaux.fr.} Partially supported by the ANR projects TRECOS ANR-20-CE40-0009 (2021--2024) and NumOpTes ANR-22-CE46-0005 (2023--2026)},  \and F. Jelassi\thanks{Alliance Sorbonne Universit\'e. UTC, EA 2222, Laboratoire de Math\'ematiques Appliqu\'ees de Compi\`egne, F-60205. E-mail: {\tt faten.jelassi@utc.fr} Supported by the ANR project  NumOpTes ANR-22-CE46-0005 (2023--2026).} \and M. Mint Brahim\thanks{CMP composites, 17, rue Jean-Baptiste Perrin 33320, EYSINES. Supported by the ANR project PHASEFIELD ANR-16--CE40-0026-01. E-mail: {\tt maimouna.brahim@cmpcomposites.fr}.}
}  

\date\today
\begin{document}
\maketitle
\begin{abstract}

We analyze an optimization problem of the conductivity in a composite material arising in a heat conduction energy storage problem. The model is described by the heat equation that specifies the heat exchange between two types of materials with  different conductive properties with Dirichlet-Neumann boundary conditions on the external part of the domain, and on the interface characterized by  the resisting coefficient between the highly conductive material and the less conductive material. The main purpose of the paper is to compute  a  shape gradient of an optimization functional in order to accurately determine  the optimal location of the conductive material using a classical   shape optimization strategy. We also present some numerical experiments to illustrate the efficiency of the proposed method. 
\end{abstract}
%
%
%
%
%
\section{Introduction}

Thermal energy storage (TES) plays a pivotal role in achieving effective and efficient heat generation and utilization, especially when there is a spatial and temporal mismatch between heat supply and demand. Among various heat storage processes, latent heat storage stands out, involving a phase transformation of storage materials known as phase change materials (PCMs), typically transitioning between solid and liquid states. Enhancing the thermo-physical properties of PCMs through special additives or composite development is recognized as a scientific challenge to bridge the gap between TES targets and current performance levels.

\medskip
Additives, particularly carrier materials, play a crucial role in improving heat and mass transfer, enhancing thermo-mechanical stability, preventing undesirable effects like segregation, or undercooling in PCMs, and improving cyclability. One common experimental approach to boost heat transfer in heat storage materials is through the addition of extended surfaces or encapsulated phase change materials. Despite the attractive energy density of PCMs, they often exhibit low thermal conductivity. Therefore, enhancing energy storage efficiency becomes imperative, and this can be achieved by incorporating spherical or other shape capsules that improve conductivity properties. Naturally, the key question revolves around studying the influence of the shape and position of capsules to accelerate the melting process.

\medskip
To answer this goal, a suitable optimization problem is considered for which different mathematical tools can be used, like geometrical and topological shape optimization based on level set method (see \cite{DelgadoPhd}).  The main idea of this Eulerian approach is to represent the shape of the composite structure in an implicit way within a fixed meshed working domain.  The level set function is then advected through Eikonal  Hamilton-Jacobi equation in which the normal velocity is  provided by the  shape derivative of the cost function.   The difficulties users have to cope with are listed by G. Allaire (see \cite{Allaire02,Allaire14} and references therein). The challenge is, usually, to perform an appropriate interpolation scheme to smoothen the discontinuous material parameters across the interfaces or to consider diffuse interfaces instead of  sharp  ones. 

\smallskip
Various techniques have been suggested in the literature, such as density-based methods (we list essentially the Solid Isotopic Material with Penalization, SIMP, topology optimization method \cite{Bendsoe02}  and  the homogenized one \cite{Tartar09}). The drawback of such methods is that they are not effective to capture precisely  interfaces unless dense girds at the vicinity of these ones are created.  The phase field method (see \cite{Zhou08}) is an interesting alternative to deal with diffuse interfaces in multimaterial structural optimization. In recent years, original approach has been developed combining level set method and Extended Finite Element  Method (XFEM) in structural topology and shape optimization field  that  further accurate an efficient issues.  G. Allaire  proposes a level set based mesh evolution method applied to multimaterials compliance minimization \cite{Allaire14}. The milestone   of this approach is to adapt mesh at each stage of the optimization process, so to obtain an exact meshed description of the shapes while benefiting the whole flexibility of the level set method. An alternative to geometrical shape optimization is the topological shape optimization, a recent mathematical technology championed by  Sokolowski \cite{SokZol92}, Masmoudi \cite{Larnier13}, etc.

\medskip
In this article, we propose to take advantage of the spherical shape of the capsule to adopt the strategy of  computing  a  shape gradient of an optimization functional in order to compute accurately the optimal location of the conductive material using a classical strategy of shape optimization. We present an investigation of the optimization of the conductivity in a composite material. We consider the  composite material as  the assembly of two or more materials that are not miscible, of a different nature, and able to combine several mechanical or chemical characteristics specific to each component. 
Doing so, the resulting material may have better overall thermo-physical properties than its constituent materials and allow for a wider range of applications (see~\cite{Sai16}). This is why in recent years, composite materials have established themselves in many cutting-edge sectors such as space, aeronautics, shipbuilding or even the automotive industry. 

\medskip
Often times the composite consists of a matrix of a given material (for example cement or carbon)
embedded with a second material to enhance a given characteristic such as thermal conductivity or 
the resistance to fractures. In such case, when subject to a heat load, a temperature drop at the interface between components of a composite is observed. We describe this discontinuity in the temperature field using a thermal resistance which is known as inter-facial thermal resistance and which is the result of two phenomena. The first is the thermal resistance observed when two components are in contact as a result of poor mechanical and chemical bounds between the components. The second is when there is a discontinuity in the thermal property of the components of the material such as the thermal conductivity~\cite{Pietrak}.

\medskip
For composites subject to high temperatures, for example in applications for thermal energy storage which is at the origin of our motivation, arrangements of fibers and their shape play a crucial role in the effectiveness of the composite. Using experimental means in order to design the position and or optimal shapes for such materials can become very expensive and time consuming. 

\medskip
It is therefore interesting to develop a theoretical argument to optimize the fiber placements within the material without  any modification of the topology, and this is precisely our goal. In Section \ref{Sec-Theoretical-Study}, we give the mathematical formalism for  the unsteady heat conduction problem in a composite media with thermal contact resistance at the interface between the components, and develop the theoretical approach for the shape optimization problem in this context. Then in Section \ref{Sec-Num-Study}, numerical examples illustrate the efficiency of the proposed method in the $2$d setting when prescribing the shape of the fiber to be a disk of prescribed radius.
\bigskip

\section{Mathematical setting and results}\label{Sec-Theoretical-Study}

In this section we will present the optimization problem we are going to analyze and then introduce the main result, the  computation of the shape gradient of the optimization functional. In addition, we consider the particular case in which the domain $\mathcal{O}$ is a ball.  We will also give the proof of the main result. 

\subsection{Statement of the problem}\label{Subsec-Setting}

We consider the domain 
$$
	\Omega = (0,1)^2
$$
and we decompose its boundary as follows:
$$
	\partial \Omega = \Gamma_0 \cup \Gamma_n, \quad \Gamma_0 = (0,1) \times \{0\}, \quad  \Gamma_n = \partial\Omega \setminus \Gamma_0.
$$

Let $\kappa >1$ be a (large) positive constant corresponding to the conductivity of the highly conductive material, and $R>0$ denotes a resisting coefficient between the highly conductive material and the less conductive material, supposed to be of conductivity one.

Given a smooth open subset $\O$ of $\Omega$ with $\O \Subset \Omega$, we set 
\begin{equation}
	\label{Def-S-from-O}
	\S = \Omega \setminus \overline{\O}
\end{equation}
and we consider the solution $u = (u_\S, u_\O)$ of the following problem: 
\begin{equation}
	\label{Eq-O}
	\left\{
		\begin{array}{ll}
			\partial_t u_\S - \Delta u_\S = 0, \quad & \hbox{ in } (0,T) \times \S, 
			\\[1mm]
			\partial_t u_\O - \kappa \Delta u_\O = 0, \quad & \hbox{ in } (0,T) \times \O, 
			\\[1mm]
			\partial_n u_\S = \kappa \partial_n u_\O, \quad & \hbox{ on } (0,T) \times \partial \O, 
			\\[1mm]
			R \partial_n u_\S = u_\O - u_\S,  \quad & \hbox{ on } (0,T) \times \partial \O, 
			\\[1mm]
			\partial_n u_\S = 0, \quad & \hbox{ on } (0,T) \times \Gamma_n, 
			\\[1mm]
			u_\S = U_M,  \quad & \hbox{ on } (0,T) \times \Gamma_0, 
			\\[1mm]
			(u_{\S} (0, \cdot), u_\O(0, \cdot)) = (0,0), \quad & \hbox{ in } \S \times \O.
		\end{array}
	\right.
\end{equation}
Here, $u_\S= u_\S(t,x)$ denotes the temperature in the less conductive material, which fills the set $\S$, $u_\O=u_\O(t,x)$ denotes the temperature in the highly conductive material, which fills the set $\O$, $n$ denotes the outward normal to $\S$, and $U_M$ denotes a positive constant temperature. 
\medskip

To analyze the well-posedness issues for \eqref{Eq-O}, it is useful to remark that the system can be fit into an abstract form as
\begin{equation}
	\label{Abstract-Eq}
	v' + \mathcal{L}_\O v = 0, \quad t \in (0,T), \qquad v(0) = - U_M, 
\end{equation}
where $v$ corresponds to $u$ by the relation $v = (v_\S, v_\O) = (u_\S  - U_M, u_\O - U_M)$, and where $\mathcal{L}_\O$ is the operator on $L^2(\S) \times L^2(\O)$ with the scalar product 
$$
	\langle (v_\S, v_\O), (w_\S, w_\O) \rangle_{L^2(\S) \times L^2(\O)}= \langle v_\S, w_\S \rangle_{L^2(\S)} +  \langle v_\O,  w_\O \rangle_{L^2(\O)},
$$

\noindent 
(which of course can be identified with $L^2(\Omega)$ with its usual scalar product) with domain 
\begin{multline}
	\label{Domaine-of-L}
	\mathscr{D} (\mathcal{L}_\O  ) = \{ (v_\S, v_\O) \in L^2(\S) \times L^2(\O), \hbox{ with } (\Delta v_\S , \Delta v_\O) \in L^2(\S) \times L^2(\O), 
	\\
	\hbox{ and } \partial_n v_\S =\kappa \partial_n v_\O \hbox{ and }  \, R \partial_n v_\S = v_\O - v_\S \hbox{ on } \partial \O, \hbox{ and } \partial_n v_\S = 0 \hbox{ on } \Gamma_n, \, v_\S = 0 \hbox{ on } \Gamma_0 \}, 
\end{multline}
defined by 
\begin{equation}
	\label{Def-L}
	\mathcal{L}_\O  (v_\S, v_\O) = (- \Delta v_\S, - \kappa \Delta v_\O).
\end{equation}

One then easily checks that $\mathcal{L}_\O$ is positive self-adjoint: for all $(v_\S, v_\O)$ and $(w_\S, w_\O)$ in $\mathscr{D} (\mathcal{L}_\O )$, 

$$
	\langle \mathcal{L}_\O (v_\S, v_\O), (w_\S, w_\O) \rangle_{L^2(\S) \times L^2(\O)}
	= 
	\int_\S \nabla v_\S \cdot \nabla w_\S \, dx 
	+
	\kappa \int_\O \nabla v_\O \cdot \nabla w_\O \, dx 
	+ 
	\frac{1}{R} \int_{\partial\O} [v] [w] \, d\sigma, 
$$
where $[v] = v_\O - v_\S$ and $[w] = w_\O - w_\S$ denote the jumps at the interface $\partial\O$.

\medskip
Accordingly, the well-posedness of \eqref{Abstract-Eq} can be derived through standard arguments, see for instance \cite[Section 7.4]{Brezis}, and we will not give more details about it here.

On the other hand, it is convenient to rewrite problem \eqref{Eq-O} in a variational form. In order to do this, we introduce the space 
$$
	H^1_\O(\Omega)  = \{ v = (v_\S, v_\O) \in L^2(\Omega) \, \text{ such that } (v_\S, v_\O) \in H^1(\S) \times H^1(\O), \quad v_s = 0 \hbox{ on } \Gamma_0 \}, 
$$
and we rewrite the problem \eqref{Eq-O} in its variational form as follows:  find $u = (u_\S,u_\O)$ such that $u- U_M \in L^2(0,T; H^1_\O(\Omega)) \cap H^1(0,T; (H^1_\O(\Omega))')$ and for all $z \in  H^1_\O(\Omega)$, we have (equality in $L^2(0,T)$)

\begin{equation}
	\label{Variational-Form-1}
	\frac{d}{dt} \left( \langle (u_\S, u_\O), (z_\S, z_\O) \rangle_{L^2(\S) \times L^2(\O)}
 \right)
 	+ 
	\int_\S \nabla u_\S \cdot \nabla z_\S \, dx 
	+
	\kappa \int_\O \nabla u_\O \cdot \nabla z_\O \, dx 
	+ 
	\frac{1}{R} \int_{\partial\O} [u] [z] \, d\sigma = 0
\end{equation}
with
\begin{equation}
	\label{Variational-Form-2}
	 	 (u_\S(0), u_\O(0)) =( 0,0).
\end{equation}
\smallskip

It is also observed, since the unique non-trivial datum is $U_M$, which is a positive constant temperature, the solution $u = (u_\S, u_\O)$ of \eqref{Eq-O} will be smooth in the following sense: $u_\S \in \mathscr{C}^\infty(\overline\S)$, $u_\O \in \mathscr{C}^\infty(\overline\O)$. This can be checked using classical PDE arguments, for instance by constructing a smooth lifting $\tilde u_M$ of $U_M$, compactly supported  in a neighborhood of $(0,T) \times \Gamma_0$ not intersecting $\O$, thus yielding $u - \tilde u_M$ as the solution of \eqref{Eq-O} with homogeneous Dirichlet boundary conditions on $(0,T) \times \Gamma_0$ and a smooth source term localized in $\S$ and away from the interface $\partial \O$. Classical strategies then apply immediately, see for instance \cite[Section 7.1.3]{EvansPDE}. 

\medskip

The optimization problem we would like to address consists in studying the following functional: 
\begin{equation}\label{Jfitness}
	J(\O) = \int_0^T \!\!\!\int_\S | u_\S(t,x)  - U_M|^2 \, dt dx, 
\end{equation}
where $u_\S$ satisfies \eqref{Eq-O} is depending on $\O$ through the definition of the sets $\mathcal{O}$ and $\mathcal{S} = \Omega \setminus \overline{\mathcal{O}}$  (see \eqref{Def-S-from-O}) and $U_M$ is the positive constant temperature appearing in the boundary condition on $\Gamma_0$ in (\ref{Eq-O}). In some sense, this functional can be seen as an evaluation of the time needed for solutions of \eqref{Eq-O} to reach the equilibrium $(u_\S, u_\O) = (U_M, U_M)$. 

Therefore, the optimization problem we would like to establish is the following: find $\O \in {\mathcal{U}}_{ad}$ such that 
\begin{equation}\label{JoptimPB}
	J(\O) \leq J(\O') \quad \forall\,  \O'\in {\mathcal U}_{ad},
\end{equation}
where $ {\mathcal U}_{ad}$ is suitable admissible set of domains for $\mathcal{O}$.


\begin{remark}
	Before going further, let us note that there could be many different functionals which may be interpreted as an evaluation of the time needed for solutions of \eqref{Eq-O} to reach the equilibrium $(u_\S, u_\O) = (U_M, U_M)$.
	
	For instance, since $\mathcal{L}_\O$ is positive self-adjoint with compact resolvent, its spectrum is composed of a sequence of eigenvalues $\lambda_1(\O) \leq \cdots \leq \lambda_j(\O) \leq \cdots \to \infty$ and corresponding eigenvectors $(\Phi_j(\O))_{j \in \N}$, which form an orthonormal basis of $L^2(\S) \times L^2(\O)$. One easily checks that for all $t\geq 0$ we have
	$$
		\| e^{- t \mathcal{L}_\O}\|_{\mathscr{L}(L^2(\S) \times L^2(\O))} = e^{-t \lambda_1(\O)}.
	$$
	Therefore, another optimization problem of interest could be the minimization of the first eigenvalue $\lambda_1(\O)$ with respect to $\O$ in a suitable class of admissible domains. We did not explore this way so far, since this is certainly a more intricate study that the one we propose here.

	Additionally, let us remark that, for the application in mind, a change of phase appears when the temperature gets close to the temperature $U_M$, which yields to a non-linear parabolic model instead of \eqref{Eq-O}. With such non-linear equations, the above interpretation of the time needed to reach the equilibrium in terms of the eigenvalue of the operator is not so evident. \hfill $\blacksquare$
\end{remark}

\subsection{Computation of the shape derivative: main result}\label{Subsec-Shape-Der}

For numerical approximation of the optimization problem \eqref{JoptimPB}, we  compute the gradient of the functional $J$ at some open set $\O_0\Subset \Omega$. Following the classical 
strategy of shape optimization, see e.g. \cite{Henrot-Pierre-2005}, we somehow embed a differential structure on the set of all open sets. This can be done as follows. 

\medskip 
Given a smooth vector field $f : \Omega \to \R^2$ such that $f$ vanishes in a neighborhood of $\partial \Omega$, we introduce $X_f$ as the flow corresponding to $f$, that is the solution of the following ordinary differential equations:
\begin{equation}\label{difpb}
\begin{cases}
	\dfrac{dX_f}{ds}(s, x) = f(X_{f}(s,x)), \quad &s \in \R,\, x \in \Omega, 
	\\[2mm]
	X_f(0, x) = x, \quad & x \in \Omega.
	\end{cases}
\end{equation}

%
%

\medskip

We then consider the family of open sets given by $s \mapsto \O_{f,s} = X_f(s, \O_0)$ and the following quantity: 
\begin{equation}
	\label{Shape-Derivative}
	\frac{d}{ds} \left( J (\O_{f,s}) \right) |_{s = 0}, 
\end{equation}
which, in some sense, corresponds to the derivative of $J$ at $\O_0$ in the direction of (the deformation given by) $f$.
\begin{remark}
	Note that another possibility could be to set $\tilde\O_{s f} = (Id + s f)(\O_0)$ for $s$ in a neighborhood of $0$ and to analyze the derivability of $s \mapsto J (\tilde \O_{sf})$ at $0$, corresponding to $\tilde X_f(s,x) = x + s f(x)$ in a neighborhood of $\partial \O_0$. Note that the two computations yield the same result for the derivative, since what matters is the first order variation of $\O_0$ under the vector field $f$, see for instance \cite[Remark 5.2.9]{Henrot-Pierre-2005}.    \hfill $\blacksquare$                        
\end{remark} 

According to the structure theorem \cite[Proposition 5.9.1]{Henrot-Pierre-2005}, we know that, if $s \mapsto J (\O_{f,s})$ is differentiable at $s = 0$ for all $f$, then we should have:
\begin{equation}
	\label{Shape-Derivative-2}
	\frac{d}{ds} \left( J (\O_{f,s}) \right) |_{s = 0} = \ell ( f\cdot n|_{\partial \O}) , 
\end{equation}
where $\ell$ is a linear form. 

Our goal is to identify this linear form. More precisely, we show the following result:
%

\begin{theorem}\label{ThmMain}
	Let $\O_0$ be a smooth open domain of $\Omega$ with $\O_0 \Subset \Omega$, and $f$ be a smooth vector field vanishing in a neighborhood of $\partial\Omega$. For convenience, we next denote $\O_0$ and $\S_0 = \Omega \setminus \O_0$ simply by $\O$ and $\S$.
	
	Let us denote by $u = (u_\S,u_\O)$ the solution of \eqref{Eq-O}, and introduce $g = (g_\S, g_\O)$ the solution of the adjoint problem
	\begin{equation}
	\label{Eq-O-Adj}
	\left\{
		\begin{array}{ll}
			-\partial_t g_\S - \Delta g_\S = (u_\S - U_M), \quad & \hbox{ in } (0,T) \times \S, 
			\\[1mm]
			- \partial_t g_\O - \kappa \Delta g_\O = 0, \quad & \hbox{ in } (0,T) \times \O, 
			\\[1mm]
			\partial_n g_\S = \kappa \partial_n g_\O, \quad & \hbox{ on } (0,T) \times \partial \O, 
			\\[1mm]
			R \partial_n g_\S = g_\O - g_\S,  \quad & \hbox{ on } (0,T) \times \partial \O, 
			\\[1mm]
			\partial_n g_\S = 0, \quad & \hbox{ on } (0,T) \times \Gamma_n, 
			\\[1mm]
			g_\S = 0,  \quad & \hbox{ on } (0,T) \times \Gamma_0, 
			\\[1mm]
			(g_{\S} (T, \cdot), g_\O(T, \cdot)) = (0,0), \quad & \hbox{ in } \S \times \O.
		\end{array}
	\right.
	\end{equation}

	We also introduce the smooth vector field $\tau : \partial \O \to \R^2$ which is such that for all $x \in \partial \O$, $\tau(x)$ is a unit tangent vector of $\partial \O$ at $x$ (the orientation is arbitrary).
	
	\smallskip 
	Then we have the following  formula:
	
	\begin{multline}\label{Grad}
	\frac{d}{ds} \left( J (\O_{f,s}) \right) |_{s = 0}
	= 
	 \int_{\partial\O} f \cdot n \int_0^T \Bigg( ( u_\S - U_M)^2
		+2\Big(\frac{\kappa - 1}{\kappa}\Big) \frac{1}{R^2}(g_\O - g_\S)(u_\O - u_\S) 
	\\
		   - 2 \tau \cdot \nabla a -2 a \divO \tau  
		 +2 \kappa D^2 u_\O n \cdot n\, g_\O -2 D^2 u_\S n\cdot n \,g_\S 
	\Bigg) \, dt \, d\sigma,
	\end{multline}
	where $a$ is given by 
	
	\begin{equation}
		\label{Def-A}
		a = g_\S \tau \cdot \nabla u_\S - \kappa g_\O  \tau \cdot \nabla u_\O \quad \hbox{ on } (0,T) \times \partial\O,  
	\end{equation}

\noindent
	and $\divO$ is the tangential divergence of the vector field $\tau$, i.e. if $\tilde \tau$ is an extension of $\tau$ in a neighborhood of $\partial\O$, $\divO \tau= \div \tilde\tau - D \tilde\tau n \cdot n = \divO \tau$.
\end{theorem}

\medskip
We will give the proof of Theorem~\ref{ThmMain} in Section~\ref{sec.proof}. 
\begin{remark}
	\label{Rem-Tangential-Div}
	 The tangential divergence of $\tau$ does not depend on the extension $\tilde \tau$ chosen to extend $\tau$ in a neighborhood of $\partial\O$ and is intrinsic (see \cite[Definition 5.4.6]{Henrot-Pierre-2005}).  \hfill $\blacksquare$
\end{remark}

%

\medskip
Let us note that in the formula \eqref{Grad}, the term 
\begin{multline*}
	\int_0^T \Bigg( ( u_\S - U_M)^2
		+2\Big(\frac{\kappa - 1}{\kappa}\Big) \frac{1}{R^2}(g_\O - g_\S)(u_\O - u_\S) 
	\\
		   - 2 \tau \cdot \nabla a -2 a \divO \tau  
		 +2 \kappa D^2 u_\O n \cdot n\, g_\O -2 D^2 u_\S n\cdot n \,g_\S 
	\Bigg) \, dt
\end{multline*}
can be interpreted as the gradient of $J$ in the following sense: if one starts from a smooth open domain $\O$ included in $\Omega$, to make the functional $J$ decays, one should choose a small perturbation of $\O$ of the form $X_f(s,\O)$ for $s$ small (or $(Id + s f) (\O)$) for $f$ chosen such that 
\begin{multline*}
	f\cdot n =  -\int_0^T \Bigg( ( u_\S - U_M)^2
		+2\Big(\frac{\kappa - 1}{\kappa}\Big) \frac{1}{R^2}(g_\O - g_\S)(u_\O - u_\S) 
		\\
		   - 2 \tau \cdot \nabla a -2 a \divO \tau  
		 +2 \kappa D^2 u_\O n \cdot n\, g_\O -2 D^2 u_\S n\cdot n \,g_\S 
	\Bigg) \, dt,  
\end{multline*}
with $a$ as in \eqref{Def-A}.

\medskip 
When a large number of degrees of freedom is allowed to parametrize the shape of $\O$, the formula \eqref{Grad}, which requires only the computation of the solution of \eqref{Eq-O} and of the solution of the adjoint equation \eqref{Eq-O-Adj}, is then much more efficient than doing a derivative free optimization of $J$, which can become very costly for high number of parameters. 

The drawback is that the derivative \eqref{Grad} requires the computation of the traces of $D^2 u_\S$ and $D^2 u_\O$ on the interface $(0,T) \times \partial\O$, which are delicate to compute and may require a highe order numerical method to be computed properly.

\bigskip
\subsection{Shape derivative when imposing that $\O$ is a ball}\label{Subsec-Computation-Ball}

Let us consider the particular case of the shape derivative given by the formula \eqref{Grad} when we impose that the set $\O$ is a ball of prescribed radius $r>0$. We show below that the expression of the shape derivative in \eqref{Grad} can be slightly simplified. We will use it in our numerical experiments presented in Section~\ref{Sec-Num-Study}.

Indeed, when $\O = B(x_0, r)$ is the ball of the center $x_0= (x_{0,1},x_{0,2})$ and radius $r$, we can choose on $\O$ as tangential and normal vector fields to $\partial\O$ the following vectors, whose formula can be extended in a neighborhood of $\partial\O$ as follows:

    \[
    \tau (x) = \frac{1}{|x -x_0|} \begin{pmatrix} x_2 - x_{0,2} 
    \\[2mm]
    - (x_1 - x_{0,1})
    \end{pmatrix}  =\frac{ (x- x_0)^\perp}{|x -x_0|}, 
    \qquad  
  	n(x) = - \frac{1}{|x -x_0|} \begin{pmatrix} x_1 - x_{0,1} 
    \\[2mm]
    x_2 - x_{0,2}
    \end{pmatrix}=  - \frac{(x-x_0)}{|x -x_0|}, 
    \]
so that 

    \begin{equation*}
    	 \div \tau = 0, \quad D \tau n \cdot n = 0, \quad \divO \tau = 0, \quad D \tau \tau =  \frac{n}{|x-x_0|}.
	\end{equation*}	

	\noindent
	Accordingly, for
	\begin{equation*}
	         a 
	         = \tau \cdot (g_\S  \nabla u_\S- \kappa g_\O \nabla u_\O), 
        \end{equation*}
	we get
	
	\begin{align*}
        \tau \cdot  \nabla a  
       & = 
       (\tau \cdot \nabla g_\S) (\tau \cdot \nabla u_\S)
       - 
       \kappa (\tau \cdot \nabla g_\O) (\tau \cdot \nabla u_\O)
       + 
       \frac{1}{r} (g_\S \partial_n u_\S - \kappa g_\O \partial_n u_\O)
       \\[1mm]
      &  \phantom{--} + 
       g_\S D^2 u_\S \tau \cdot \tau - \kappa g_\O D^2 u_\O \tau \cdot \tau
       \\[1mm]
       &  = 
         (\tau \cdot \nabla g_\S) (\tau \cdot \nabla u_\S)
       - 
       \kappa (\tau \cdot \nabla g_\O) (\tau \cdot \nabla u_\O)
       + 
       \frac{1}{r R} (g_\S- g_\O)(u_\O - u_\S)
         \\[1mm]
       &  \phantom{--}  + 
       g_\S D^2 u_\S \tau \cdot \tau - \kappa g_\O D^2 u_\O \tau \cdot \tau. 
    \end{align*}
    Since for a square matrix $A \in \R^{2\times 2}$, $A \tau \cdot \tau + A n \cdot n = \mathop{tr\,}(A)$, it follows that
	\begin{align}
	& \int_0^T \Bigg( ( u_\S - U_M)^2
		+2\Big(\frac{\kappa - 1}{\kappa}\Big) \frac{1}{R^2}(g_\O - g_\S)(u_\O - u_\S) 
		 \notag \\
	& \qquad  \qquad  \qquad  - 2 \tau \cdot \nabla a -2 a \divO \tau  
		 +2 \kappa D^2 u_\O n \cdot n\, g_\O -2 D^2 u_\S n\cdot n \,g_\S \Bigg)
		 \, dt \notag
	\\
	&= \int_0^T \Bigg( ( u_\S - U_M)^2
		+2\Big(\frac{\kappa - 1}{\kappa R^2}+ \frac{1}{rR}\Big) (g_\O - g_\S)(u_\O - u_\S) 
		\label{Grad-ball}\\
	& \qquad  \qquad  \qquad + 2  \kappa (\tau \cdot \nabla g_\O) (\tau \cdot \nabla u_\O) - 2    (\tau \cdot \nabla g_\S) (\tau \cdot \nabla u_\S)
		 +2 \kappa \Delta u_\O g_\O -2 \Delta u_\S \,g_\S \Bigg)
		 \, dt. 	 
		 \notag
    	\end{align}
	Clearly, to compute it efficiently, even for balls, this requires the computation of traces of Laplacian for solutions $(u_\O, u_\S)$ and thus a refined precision close to the interface $\partial\O$.

\bigskip
\bigskip

\subsection{Proof of Theorem \ref{ThmMain}}\label{sec.proof}

The computation of the derivative \eqref{Shape-Derivative} is done in several steps, that we briefly explain below.
\medskip 

For convenience, we fix a smooth vector field $f \in W^{2,\infty}(\Omega)$ vanishing in a neighborhood of $\partial\Omega$. To do the computations, we express all the quantities $J(X_f(s,\O))$, where $X_f(s, \cdot)$ is the diffeomorphism in \eqref{difpb}, in the fixed domain $\O$. For simplifying the notation, we will also simply denote $X_f$ by $X$. 
We thus express $J(X(s,\O))$ in terms of 
$$
	v^s (t,y) = u^s( t, X(s,y)), \quad y \in \O, \quad \text{ or equivalently } \quad u^s (t,x) = v^s(t, X(-s,x)), \quad x \in X(s,\O),
$$
where $u^s$ denotes the solution of \eqref{Eq-O} with $\O_{s} = X(s, \O)$. This gives
\begin{align}
	J(X(s,\O)) 
	&= \int_0^T\!\!\! \int_{X(s, \S)} | u^s(t, x) - U_M|^2 \, dt dx
		\\[1mm]
	& = \int_0^T\!\!\!  \int_{\S} |u^s(t, X(s, y)) - U_M|^2 \, |\det DX(s,y)| \, dt \,dy 
		\\[1mm]
 	& = \int_0^T\!\!\!  \int_{\S} |v^s(t, y) - U_M|^2 \, |\det DX(s,y)| \, dt \,dy.
	\label{Rewriting-J}
\end{align}
Accordingly, it will be important to compute $v^s$ and its derivative at 
$s= 0$, which will be the next steps.

\bigskip
\noindent
{\bf Step 1:  Approach 1. Direct computation of $v^s$ from the variational formulation (\ref{Variational-Form-1})--(\ref{Variational-Form-2}).} 

\medskip
We start by remarking that for all $s$ small enough, $(u_\S^s, u_\O^s)$ satisfies: for all $z^s \in H^1_{\O_s}(\Omega)$, we have

\begin{multline}
	\label{Variational-Form-1-s}
	\frac{d}{dt} \left( \langle (u_\S^s, u_\O^s), (z_\S^s, z_\O^s) \rangle_{L^2(\S_s) \times L^2(\O_s)}
 \right)
 	\\[1mm]
	+ 
	\int_{\S_s} \nabla u_\S^s \cdot \nabla z_\S^s \, dx 
	+
	\kappa\int_{\O_s} \nabla u_\O^s \cdot \nabla z_\O^s \, dx 
	+ 
	\frac{1}{R} \int_{\partial\O_s} [u^s] [z^s] \, d\sigma_s = 0. 
\end{multline}
Accordingly, for $z \in H^1_{\O}(\Omega)$, since $z^s(x) = z(X(-s,x))$ belongs to $H^1_{\O_s}(\Omega)$, the change of variable formula gives (still in $L^2(0,T)$)
\begin{multline}
	\label{Variational-Form-1-s-bis}
	\frac{d}{dt} \left( \int_\Omega v^s(t,y) z (y) J_s(y)\, dy \right)
 	\\
	+ 
	\int_{\S} A_s(y) \nabla v_\S^s \cdot \nabla z_\S \, dy 
	+
	\kappa \int_{\O}  A_s(y) \nabla v_\O^s \cdot \nabla z_\O \, dy 
	+ 
	\frac{1}{R} \int_{\partial\O} [v^s] [z] \, J^s(y) d\sigma = 0, 
\end{multline}
where 

\begin{equation}
	\label{Various-Quantities}
	\left\{\begin{array}{llrl}
	J_s(y) &= \det DX(s,y), \quad &&y \in \Omega,
	\\[2mm]
	A_s(y) &= J_s(y) DX(s,y)^{-1} (DX(s,y)^t)^{-1}, \quad &&y \in \Omega,
	\\[2mm]
	J^s(y) & = J_s (y) | (DX(s,y)^t)^{-1} n_y |, \quad &&y \in \partial\O,
	\end{array}\right.
\end{equation}
 see \cite[Proposition 5.4.3]{Henrot-Pierre-2005} for the last formula.

\medskip
This variational formulation can be interpreted as follows: 

\begin{equation}
	\label{Eq-O-vs}
	\left\{
		\begin{array}{ll}
			J_s(y) \partial_t v_S^s - \div (A_s(y) \nabla v_\S^s) = 0, \quad & \hbox{ in } (0,T) \times \S, 
			\\[2mm]
			J_s(y)\partial_t v_\O^s - \kappa \div (A_s(y) \nabla v_\O^s)= 0, \quad & \hbox{ in } (0,T) \times \O, 
			\\[2mm]
			A_s (y) \nabla v^s_\S \cdot n 
			= 
			\kappa A_s (y) \nabla v^s_\O \cdot n , \quad & \hbox{ on } (0,T) \times \partial \O, 
			\\[2mm]
			R A_s (y) \nabla v^s_\S \cdot n  = (v_\O^s - v_\S^s)J^s(y),  \quad & \hbox{ on } (0,T) \times \partial \O, 
			\\[2mm]
			\partial_n v_\S^s = 0, \quad & \hbox{ on } (0,T) \times \Gamma_n, 
			\\[2mm]
			v_\S^s = U_M,  \quad & \hbox{ on } (0,T) \times \Gamma_0, 
			\\[2mm]
			(v_{\S}^s (0, \cdot), v_\O^s(0, \cdot)) = (0,0), \quad & \hbox{ in } \S \times \O.
		\end{array}
	\right.
\end{equation}

\bigskip\noindent
{\bf  Step 2: Approach 2. Direct computations of $v^s$ from the equation (\ref{Eq-O})}. 

\medskip
We have that 
\begin{align*}
	& \nabla_x u^s(t,x)  = D_x X(-s,x)^t \nabla_y v^s(t,X(-s,x))
	= 
	\left(\sum_k \partial_{y_k} v^s (t,X(-s,x)) \partial_{x_j} X_k(-s,x)\right)_j. 
\end{align*}

\noindent 
Then, we can write 

\begin{align*}
	& \Delta_x u^s(t,x) = 	\tr( D_x X(-s,x)^t D^2_y v^s(t,X(-s,x)) D_x X(-s,x)) + \Delta_x X(-s,x)\cdot \nabla_y v^s (t,X(-s,x))
	\\[1mm]
	&\hspace{0.4cm} = \sum_{j,k,\ell} \partial_{y_ky_\ell} v^s(t,X(-s,x)) \partial_{x_j} X_k(-s,x)\partial_{x_j} X_\ell(-s,x) + \sum_{j,k} \partial_{y_k} v^s(t,X(-s,x)) \partial_{x_j x_j} X_k(-s,x)
\end{align*}

\noindent 
and

\[
n_x  = \frac{D_x X(-s,X(s,y))^t n_y}{| D_x X(-s,X(s,y))^t n_y|} \quad \text{ with } y = X(-s,x) \text{ for } x \in \partial \O_s,
\]
where, for a matrix $A$, $A^t$ denotes its transpose. Note that the last identity is standard and can be deduced for instance from \cite[Proposition 5.4.15]{Henrot-Pierre-2005}.

\medskip
Accordingly, the equations on $u^s$ read as follows on $v^s =(v^s_{\S}, v^s_{\O})$:
\begin{equation}
	\label{Eq-v}
	\left\{
		\begin{array}{ll}
			\!\partial_t v^{s}_\S -\! \Big(
			\tr( D_x X(-s,X(s,y))^t D^2_y v^s_\S(t,y) D_x X(-s,X(s,y))) + \Delta_x X(-s,X(s,y))\cdot \nabla_y v^s_\S (t,y)\Big)= 0&
			\\
&\hspace{-3cm} \hbox{ in } (0,T) \times \S, 
	\\[2mm]
			\!\partial_t v^s_{\O} - \!\!\kappa 
			\Big(
			\tr( D_x X(-s,X(s,y))^t D^2_y v^s_\O(t,y) D_x X(-s,X(s,y))) + \Delta_x X(-s,X(s,y))\cdot \nabla_y v^s_\O (t,y)\Big) = 0 &
			\\
& \hspace{-3cm} \hbox{ in } (0,T) \times \O, 
\\[2mm]
\displaystyle
D_x X(-s,X(s,y))^t n_y \cdot D_x X(-s,X(s,y))^t \nabla_y v^s_\S(t,y)
\\[2mm]
\displaystyle
\hspace{3cm} 
= 
\kappa D_x X(-s,X(s,y))^t n_y\cdot D_x X(-s,X(s,y))^t \nabla_y v^s_\O(t,y)
&\hspace{-3cm} \hbox{ on } (0,T) \times \partial \O, 
\\[3mm]
\displaystyle
\frac{D_x X(-s,X(s,y))^t n_y}{| D_x X(-s,X(s,y))^t n_y|}  \cdot D_x X(-s,X(s,y))^t \nabla_y v^s_\S(t,y)
= \dfrac{1}{R} (v^s_{\O} - v^s_{\S})
&\hspace{-3cm}\hbox{ on } (0,T) \times \partial \O, 
\\[5mm]
	\partial_n v^s_{\S} = 0, \quad & \hspace{-3cm} \hbox{ on } (0,T) \times \Gamma_n, 
	\\[3mm]
v^s_{\S} = U_M,  \quad & \hspace{-3cm} \hbox{ on } (0,T) \times \Gamma_0, 
\\[3mm]
(v^s_{\S} (0, \cdot), v^s_{\O}(0, \cdot)) = (0,0), \quad &\hspace{-3cm} \hbox{ in } \S \times \O.
		\end{array}
	\right.
\end{equation}		
\begin{remark}
	It is not clear at first glance that both equations \eqref{Eq-O-vs} and \eqref{Eq-v} coincide, but this is in fact the case. Indeed, writing $\tilde A_s(y) = DX(s,y)^{-1}(DX(s,y)^t)^{-1}$, and remarking that it coincides with the matrix $DX(-s,X(s,y)) DX(-s,X(s,y))^t$, tedious computations yield that for all $j \in \{1, 2\}$, 
$$
	\sum_{i = 1}^2 
	\left(\partial_i (\tilde A_s(y))_{i,j} + \frac{\partial_i J_s(y)}{J_s(y)} (\tilde A_s(y))_{i,j} \right) 
	=
	\Delta_x X_j(-s, X(s,y)).
$$
From this identity, writing the space operator in \eqref{Eq-O-vs} under the form $\frac{1}{J_s} \div ( J_s \tilde A_s \nabla \cdot )$, it is easily seen that both equations \eqref{Eq-O-vs} and \eqref{Eq-v} are identical. \hfill $\blacksquare$
\end{remark}

\bigskip\noindent
{\bf Step 3: Derivation of $v^s$ at $s = 0$.} 

\medskip
Under the form \eqref{Eq-O-vs}, it is easy to adapt the classical argument of implicit function theorem to prove that $v^s$ is $C^1$ in a neighborhood of $s = 0$ in $L^2(0,T; H_\O^1(\Omega)) \cap H^1(0,T;( H_\O^1(\Omega))')$, see for instance \cite[Proof of Theorem 5.3.2]{Henrot-Pierre-2005} for a similar argument. 

We then compute its derivative $w = (w_\S,w_\O)$ at $s = 0$. In order to do that, we simply perform a Taylor expansion of the quantities appearing in \eqref{Various-Quantities}: 

\begin{align*}
	& J_s(y) = 1+ s \div f + O(s^2),
	\\[2mm]
	& A_s(y) =  J_s(y) (I - s (Df + Df^t) + O(s^2)), 
	\\[2mm]
	&
	J^s(y) = J_s(y) (1 - s Df n \cdot n + O(s^2)),
\end{align*}
where $O(s^2)$ are functions which are bounded by $Cs^2$ as $s \to 0$.

Accordingly, $w$ solves the variational formulation (always in $L^2(0,T)$): for all $z \in H_\O^1(\Omega)$, we have

\begin{equation}
\begin{split}
	\label{Variational-Form-w-0}
	\frac{d}{dt} \left( \int_\Omega w z  + u \div(f) z \, dy \right)
	 & + 
	\int_{\S}  (\nabla w_\S \cdot \nabla z_\S +(\div(f) \nabla u_\S - (Df + Df^t) \nabla u_\S) \cdot \nabla z_\S )\, dy 
	\\[2mm]
	& +
	\kappa \int_{\O} (\nabla w_\O \cdot \nabla z_\O +(\div(f) \nabla u_\O - (Df + Df^t) \nabla u_\O) \cdot \nabla z_\O ) \, dy 
	\\[2mm]
	& + 
	\frac{1}{R} \int_{\partial\O} ([w] + (\div(f) - Df n\cdot n)[u]) [z] \, d\sigma = 0. 
\end{split}
\end{equation}
Thus, subtracting the variational formulation \eqref{Variational-Form-1} satisfied by $u$ and applied to $\div(f) z$, we get, for all $z \in H_\O^1(\Omega)$,

\begin{multline}
	\label{Variational-Form-w}
	\frac{d}{dt} \left( \int_\Omega w z  \, dy \right)
	+ 
	\int_{\S}  (\nabla w_\S \cdot \nabla z_\S) \, dy 
	-
	 \int_{\S}( \nabla u_\S \nabla \div( f)  z_\S + (Df + Df^t) \nabla u_\S \cdot \nabla z_\S)\, dy 
	\\
	+
	\kappa \int_{\O}  (\nabla w_\O \cdot \nabla z_\O) \, dy 
	-
	\kappa \int_{\O}( \nabla u_\O \nabla \div( f)  z_\O + (Df + Df^t) \nabla u_\O \cdot \nabla z_\O)\, dy 
	\\
	+ 
	\frac{1}{R} \int_{\partial\O} [w] [z] \, d\sigma 
	-
	\frac{1}{R} \int_{\partial\O} (Df n\cdot n)[u] [z] \, d\sigma = 0. 
\end{multline}

\noindent
It is then easy to check that $w$ satisfies the following:

\begin{equation}
	\label{Eq-W}
	\left\{
		\begin{array}{ll}
			\partial_t w_\S - \Delta w_\S + 2 D^2 u_\S : Df + \nabla u_\S \cdot \Delta f = 0, \quad & \hbox{ in } (0,T) \times \S, 
			\\[2mm]
			\partial_t w_\O - \kappa \Delta w_\O + 2 \kappa D^2 u_\O : Df + \kappa \nabla u_\O \cdot \Delta f = 0, \quad & \hbox{ in } (0,T) \times \O, 
			\\[2mm]
			\partial_n w_\S -  (Df+Df^t) \nabla u_\S \cdot n = \kappa \partial_n w_\O - \kappa (Df+Df^t) \nabla u_\O \cdot n, \quad & \hbox{ on } (0,T) \times \partial \O, 
			\\[2mm]
			\ds
			\partial_n w_\S -  (Df+Df^t) \nabla u_\S \cdot n
			= \frac{1}{R} (w_\O - w_\S) -\frac{Df n\cdot n}{R} (u_\O - u_\S),  \quad & \hbox{ on } (0,T) \times \partial \O, 
			\\[2mm]
			\partial_n w_\S = 0, \quad & \hbox{ on } (0,T) \times \Gamma_n, 
			\\[2mm]
			w_\S =0,  \quad & \hbox{ on } (0,T) \times \Gamma_0, 
			\\[2mm]
			(w_{\S} (0, \cdot), w_\O(0, \cdot)) = (0,0), \quad & \hbox{ in } \S \times \O. 
		\end{array}
	\right.
\end{equation}		

\noindent 
Note that, denoting by $\tau$ a tangential vector field to $\partial\O$, using the interface conditions satisfied by $u$ on $\partial\O$, the interface conditions in \eqref{Eq-W}  can be rewritten as 

\begin{equation}
	\label{Eq-W-Interface}
	\left\{
		\begin{array}{ll}
			\partial_n w_\S -  ((Df+Df^t) \tau \cdot n) \tau\cdot \nabla u_\S  = \kappa \partial_n w_\O - \kappa ((Df+Df^t) \tau \cdot n) \tau\cdot \nabla u_\O, & \hbox{ on } (0,T) \times \partial \O, 
			\\[2mm]
			\ds \partial_n w_\S - ((Df +Df^t)\tau \cdot n) \tau\cdot \nabla u_\S - \partial_n u_\S (Df n \cdot n)= \frac{1}{R} (w_\O - w_\S), & \hbox{ on } (0,T) \times \partial \O. 
		\end{array}
	\right.
\end{equation}	

\bigskip \noindent 
{\bf Step 4: Computing the derivative of $s \mapsto J(X(s,\O))$ at $s = 0$.} 

\medskip
Recalling \eqref{Rewriting-J}, we have

\begin{align}
	\frac{d}{ds} \left( J(X(s,\O)) \right)
 	& = 
	2 \int_0^T\!\!\!  \int_{\S} w_\S (u_\S - U_M)  \, dt \,dy
	+
	\int_0^T\!\!\!  \int_{\S} |u_\S - U_M|^2 \, \div(f) \, dt \,dy
	\nonumber\\[2mm]
	&
	= 
	2 \int_0^T\!\!\!  \int_{\S} (w_\S - f \cdot \nabla u_\S) (u_\S - U_M)  \, dt \,dy
	+
	\int_0^T\!\!\!  \int_{\S}  \div\left(f(u_\S - U_M)^2\right) \, dt \,dy
       \nonumber \\[2mm]
	&
	= 
	2 \int_0^T\!\!\!  \int_{\S} (w_\S  -f \cdot \nabla u_\S) (u_\S - U_M)  \, dt \,dy
	+
	\int_0^T\!\!\!  \int_{\partial\O} f\cdot n (u_\S - U_M)^2 \, dt \,d\sigma.
	\label{Derivative-J-s=0}
\end{align}
Using the solution $g = (g_\S, g_\O)$ of the adjoint equation \eqref{Eq-O-Adj}, the first term in the derivative of $s \mapsto J(X(s,\O))$ at $s = 0$ can be rewritten as follows: 
\begin{align}
	& \int_0^T\!\!\!  \int_{\S} (u_\S - U_M) ( w_\S - f \cdot \nabla u_\S) \, dt \, dy
	\notag
	\\[2mm]
	& = 
	 \int_0^T\!\!\!  \int_\S ( - \partial_t g_\S - \Delta g_\S)  ( w_\S - f \cdot \nabla u_\S) \, dt dy
	 + 
	 \int_0^T\!\!\!  \int_\O ( - \partial_t g_\O - \kappa \Delta g_\O)  ( w_\O - f \cdot \nabla u_\O) \, dt dy 
	\notag
	 \\[2mm]
	 & 
	 =  \int_0^T \!\!\! \int_\S g_\S ((\partial_t - \Delta) (w_\S - f \cdot \nabla u_\S)) \, dt dy
	 + 
	 \int_0^T\!\!\!  \int_\O g_\O ((\partial_t - \kappa \Delta) (w_\O - f \cdot \nabla u_\O)) \, dt dy
	 \notag
	 \\[2mm]
	 &
	 \phantom{--}  - \int_0^T\!\!\!  \int_{\partial O} 
	\left( 
	( \partial_n g_\S (w_\S - f \cdot \nabla u_\S)-  \kappa \partial_n g_\O (w_\O - f \cdot \nabla u_\O) ) 
	\right)\, dt d\sigma
	\label{Derivative-With-Adjoints-0}
	\\[2mm]
	& 
	 \phantom{--} - \int_0^T \!\!\! \int_{\partial O} 
	\left( 
	(\kappa \partial_n (w_\O  - f \cdot \nabla u_\O) g_\O- \partial_n (w_\S -f \cdot \nabla u_\S) g_\S)
	 	\right)\, dt d\sigma. 
		\notag
\end{align}
We next compute explicitly each term and simplify the terms in the right hand side of the identity in \eqref{Derivative-With-Adjoints-0} as much as possible.

First, we observe that  \eqref{Eq-W} implies that
\begin{equation}
	\label{Internal-Terms}
	(\partial_t - \Delta) (w_\S - f \cdot \nabla u_\S) = 0 \quad  \text{ in } (0,T) \times \S, 
	\qquad
	 (\partial_t - \kappa \Delta) (w_\O - f \cdot \nabla u_\O) = 0  \quad  \text{ in } (0,T) \times \O.
\end{equation}
We should thus focus on the interface terms (IT in the following). Using the boundary conditions $\partial_n g_\S = \kappa \partial_n g_\O = (g_\O - g_\S)/R$ on the interface $(0,T) \times \partial\O$, we get:

\begin{align*}
    IT  = &  - ( \partial_n g_\S (w_\S - f \cdot \nabla u_\S)+ \kappa \partial_n g_\O (w_\O - f \cdot \nabla u_\O) ) 
	\\[2mm]
	  & \quad -  (\kappa \partial_n (w_\O  - f \cdot \nabla u_\O) g_\O + \partial_n (w_\S -f \cdot \nabla u_\S) g_\S)
	\\[2mm]
	 = &  \frac{1}{R}(g_\O - g_\S) \left( (w_\O - f \cdot \nabla u_\O) - (w_\S - f \cdot \nabla u_\S) \right)
		\\[2mm]
	  & \quad 
	- \kappa \partial_n w_\O  + \partial_n w_\S g_\S + \kappa \partial_n( f \cdot \nabla u_\O) g_\O - \partial_n (f \cdot \nabla u_\S) g_\S.	
\end{align*}
On one hand, recalling that 
$$
\partial_n w_\S -  (Df+Df^t) \nabla u_\S \cdot n 
= \kappa \partial_n w_\O - \kappa (Df+Df^t) \nabla u_\O \cdot n	
= \frac{1}{R} (w_\O - w_\S) -\frac{Df n\cdot n}{R} (u_\O - u_\S),
$$
we have
\begin{multline*}
	- \kappa \partial_n w_\O g_\O+ \partial_n w_\S g_\S
	\\
	= 
	- \frac{1}{R} (g_\O - g_\S) (w_\O - w_\S) + \frac{Df n \cdot n}{R} (g_\O - g_\S) (u_\O - u_\S)
	- \kappa (Df + Df^t) \nabla u_\O \cdot n\, g_\O +(Df + Df^t) \nabla u_\S \cdot n \,g_\S.
\end{multline*}
On the other hand, computations also show that 
\begin{multline*}
	 \kappa \partial_n (f \cdot \nabla u_\O) g_\O - \partial_n (f \cdot \nabla u_\S) g_\S 
	\\[2mm]
	= 
	\kappa Df^t \nabla u_\O \cdot n\, g_\O - Df^t \nabla u_\S \cdot n\, g_\S
	+ 
	\kappa D^2 u_\O f \cdot n \,g_\O - D^2 u_\S f \cdot n \, g_\S.
\end{multline*}

\noindent
Accordingly, we obtain
\begin{align*}
	IT &  
	= - \frac{1}{R}(g_\O - g_\S) \left( f \cdot \nabla u_\O - f \cdot \nabla u_\S) \right)
	+ \frac{Df n \cdot n}{R} (g_\O - g_\S) (u_\O - u_\S)
	\\[2mm]
	&+ Df \nabla u_\S \cdot n \,g_\S - \kappa Df \nabla u_\O \cdot n\, g_\O
	+\kappa D^2 u_\O f \cdot n \,g_\O - D^2 u_\S f \cdot n \,g_\S.
\end{align*}
Introducing a tangential vector field $\tau$ to $\partial\O$, and writing $\nabla u = n \partial_n u + \tau \partial_\tau u$, we get 
\begin{align*}
	Df \nabla u_\S \cdot n \,g_\S - \kappa Df \nabla u_\O \cdot n\, g_\O
	& = 
	Df n \cdot n  ( \partial_n u_\S g_\S - \kappa \partial_n u_\O g_\O) 
	+ 
	Df \tau \cdot n ( \partial_\tau u_\S g_\S - \kappa \partial_\tau u_\O g_\O)
	\\[2mm]
	& = 
	- \frac{Df n \cdot n}{R} (u_\O - u_\S) (g_\O - g_\S)
	+ 
	Df \tau \cdot n ( \partial_\tau u_\S g_\S - \kappa \partial_\tau u_\O g_\O),  
\end{align*}
so that 
\begin{multline}
	IT = - \frac{1}{R}(g_\O - g_\S) \left( f \cdot \nabla u_\O - f \cdot \nabla u_\S \right)
	+ Df \tau \cdot n ( \partial_\tau u_\S g_\S - \kappa \partial_\tau u_\O g_\O)
	\\
	+ \kappa D^2 u_\O f \cdot n \,g_\O - D^2 u_\S f \cdot n \,g_\S.
\end{multline}

Using the convention of implicit summation, and setting $a$ as in \eqref{Def-A}, we get from \cite[Theorem 5.4.13]{Henrot-Pierre-2005} (here, $H$ denotes the mean curvature of $\partial \O$, and we choose an extension $\tilde \tau$ of $\tau$ in a neighborhood of $\partial\O$ that we still denote $\tau$ for simplicity) that 
\begin{align*}
	\int_{\partial \O} Df \tau \cdot n \, a \, d\sigma
	& = \int_{\partial \O} \partial_j f_i \tau_j  n_i a \, d\sigma
	\\
	& = - \int_{\partial \O} f_i \partial_j( \tau_j n_i a )\, d\sigma+ \int_{\partial \O} \left(\partial_n (f_i \tau_j n_i a ) + H f_i \tau_j n_i a \right) n_j\, d\sigma
	\\ 
	& = \int_{\partial \O} f \cdot n ( - \partial_j ( \tau_j a) + D\tau n \cdot n\, a) \, d\sigma
	- \int_{\partial \O} Dn \tau \cdot f a\, d\sigma
	\\
	& = \int_{\partial \O} f \cdot n ( -\tau \cdot \nabla a - a ( \div \tau - D \tau n \cdot n)) \, d\sigma
	- \int_{\partial \O} Dn \tau \cdot f \,a\, d\sigma.
\end{align*}
Note that $\div \tau - D \tau n \cdot n = \divO \tau$ is the tangential divergence of $\tau$ and is independent of the extension chosen for $\tau$ (recall Remark \ref{Rem-Tangential-Div} and \cite[Definition 5.4.6]{Henrot-Pierre-2005}).

From this identity, we deduce that 
\begin{align*}
	\int_{\partial \O} IT \, d\sigma 
	 & =  
	 \int_{\partial \O}
	  \Big(
	\frac{1}{R} (g_\O - g_\S) f \cdot ( \nabla u_\S -   \nabla u_\O) 
	+
	   f \cdot n (  -\tau \cdot \nabla a - a  \divO \tau ) 
	- Dn \tau \cdot f \,a
	\\[1mm]
	& \qquad  +\kappa D^2 u_\O f\cdot n \,g_\O - D^2 u_\S f \cdot n \,g_\S 
	\Big)d\sigma
	\\[1mm]
	& = 
	 \int_{\partial \O}\Big(
	  f \cdot n 
	   \Big(
	    \frac{1}{R^2}(g_\O - g_\S)(u_\O - u_\S) (1 - 1/\kappa)
		   -\tau \cdot \nabla a - a  \divO \tau
		-  Dn \tau \cdot n\, a 
		\\[1mm]
	& \qquad \quad
		+ \kappa D^2 u_\O n \cdot n \,g_\O - D^2 u_\S n\cdot n \,g_\S 
		\Big) \Big)d\sigma
	\\[1mm]
	& 
	\quad + 
	 \int_{\partial \O}
	  f \cdot \tau 
	   \left(
	   \frac{1}{R} (g_\O - g_\S) \tau \cdot (  \nabla u_\S -   \nabla u_\O) 
		-  Dn \tau \cdot \tau \,a
		+ \kappa D^2 u_\O \tau \cdot n \,g_\O - D^2 u_\S \tau \cdot n\, g_\S 
		\right)d\sigma.
\end{align*}
Before going further, let us also immediately point out that %
$$ 
	Dn \tau \cdot n = \partial_j n_i \tau_j n_i = \tau \cdot \nabla \left(\frac{|n|^2}2\right) = 0.
$$

It is then interesting to consider the term involving $f\cdot \tau$, as this term should vanish for structural reasons. 

 In order to do so, we   differentiate the identity 
$$
				\partial_n u_\S = \kappa \partial_n u_\O= \frac{1}{R} (u_\O - u_\S),  \quad  \hbox{ on } (0,T) \times \partial \O, 
$$
in the direction of $\tau$, yielding the following identity
$$
	 Dn \tau \cdot \nabla u_\S + D^2 u_\S \tau \cdot n 
	= 
	\kappa (Dn \tau \cdot \nabla u_\O + D^2 u_\O \tau \cdot n ) 
	= 
	 \frac{1}{R} \tau \cdot (\nabla u_\O - \nabla u_\S).
$$

Therefore, recalling the definition of $a$ in \eqref{Def-A}, we have
\begin{align*}
	 &\frac{1}{R} (g_\O - g_\S) \tau \cdot (  \nabla u_\S -   \nabla u_\O) 
		-  Dn \tau \cdot \tau \,a
		+ \kappa D^2 u_\O \tau \cdot n \,g_\O - D^2 u_\S \tau \cdot n \,g_\S 
	\\[1mm]
	&= 
	g_\O \left( \frac{1}{R}\tau \cdot (  \nabla u_\S -   \nabla u_\O) + (Dn \tau \cdot \tau) \kappa \tau \cdot \nabla u_\O  + \kappa D^2u_\O \tau \cdot n \right)
		\\[1mm]
	& \hspace{1cm}+
	g_\S \left( \frac{1}{R}\tau \cdot (  \nabla u_\O -   \nabla u_\S) - (Dn \tau \cdot \tau)\tau \cdot \nabla u_\S - D^2 u_\S \tau \cdot n \right)
		\\[1mm]
	&= 
	\kappa g_\O (- Dn \tau \cdot n) \partial_n u_\O + g_\S (Dn \tau \cdot  n) \partial_n u_\S = 0.
\end{align*}

Accordingly the expression of the interface term simplifies and we get:
 \begin{align*}
 	\int_{\partial \O} IT \, d\sigma 
	 = 
	 \int_{\partial \O}
	  f \cdot n 
	   \Big(
	    \frac{1}{R^2}(g_\O - g_\S)(u_\O - u_\S) (1 - 1/\kappa)
		  & -\tau \cdot \nabla a - a  \divO \tau
	\\
		&+ \kappa D^2 u_\O n \cdot n \,g_\O - D^2 u_\S n\cdot n \,g_\S 
		 \Big)d\sigma.
\end{align*}

Combining \eqref{Derivative-J-s=0}, \eqref{Derivative-With-Adjoints-0}, \eqref{Internal-Terms} and the above identity, we conclude 
\begin{multline*}
	\frac{d}{ds} \left( J (\O_{f,s}) \right) |_{s = 0}
	= 
	\int_0^T\!\!\! \int_{\partial\O} f \cdot n  
	\Bigg(
	( u_\S - U_M)^2\,  
	 + 2
	   \Big(
	    \frac{1}{R^2}(g_\O - g_\S)(u_\O - u_\S) (1 - 1/\kappa)
	  \\
		   -\tau \cdot \nabla a - a  \divO \tau
		+ \kappa D^2 u_\O n \cdot n\, g_\O - D^2 u_\S n\cdot n \,g_\S 
		 \Big)\Bigg)dt\, d\sigma.
\end{multline*}
as announced in the Theorem \ref{ThmMain}. 
\hfill $ \blacksquare$

\bigskip 

%
\section{Numerical results}\label{Sec-Num-Study}
In this section we test the theoretical result of Theorem  \ref{ThmMain}. Therefore, we will consider along our numerical study  a  geometry composed of a less conductive homogeneous material 
in a unit square sample containing a disc  filled with a more conductive material that acts like a catalyst to speed up the heat conduction process. 

\medskip
A gradient algorithm based on the shape derivative adapted  for the ball case and presented in Section~\ref{Subsec-Computation-Ball} will be used to obtain an optimal geometric layout. For all our  test cases, the state problem and the adjoint problem are solved using mixed finite element method with the first order finite Raviart-Thomas finite elements for the heat flux and a piece-wise constant finite elements for the temperature field (see \cite{FBB,MMB,Raviart}). The first order Euler backward method is used for time stepping. The presented results are obtained using the free finite element software \texttt{Freefem++} (see \cite{Hecht}).

\medskip
Since we are imposing that $\O$ is a ball of prescribed radius, we will use the computations in Section~\ref{Subsec-Computation-Ball}. Also note that, to guarantee that $\O$ stays a ball at each step of the optimization process, we will choose $X_f$ under the form of translations, that is $f$ a constant vector close to $\partial\O$. Note that here, rotations are irrelevant since rotation a ball does not modify its shape. Mathematically,  this correspond to a vector field $f$ for which  $f \cdot n = 0$. 

\medskip
Accordingly, at each step of the optimization process, we will use the formula \eqref{Grad-ball}, in which we will use $\Delta u_\S =\partial_t u_\S$ and $\kappa \Delta u_\O = \partial_t u_\O$ on $(0,T) \times \partial\O$, which is true by taking the traces of the equations of $u_\S$ and $u_\O$ in  \eqref{Eq-O} on the interface $(0,T) \times \partial \O$, and which seems to give more precise results in our simulations.

\bigskip
\subsection{A validation case}\label{sec.numvc}
The purpose of the first test case is to validate the implementation of the method. For this we will consider a  fabricated desired solution that we aim to reproduce thanks to our algorithm. 
This strategy requires to consider a desired temperature state which depending on the time and space coordinates. Although this case is outside the mathematical framework of Theorem \ref{ThmMain}, the proof of Theorem \ref{ThmMain} can be adapted to this test case easily,  and this test case is good to validate our algorithm.

\medskip
More precisely, we consider a desired geometrical layout with the disc at the top of the square where the center coordinates of the disc are $x_D=0.5$ and $y_D=0.75$ and the radius is $r=0.2$. The corresponding temperature field at different time step is calculated 
by solving \eqref{Eq-O} with $\kappa = 100$, $R = 10^{-2}$, $U_M = 500$ and $T = 0.5$ duration of the simulation. This information is stored as $u_D(t,x)$ and passed on into the optimization problem: Minimize, among the sets $\O$ which are disks of radius $r$, the functional 
$$
    J(\O) = \int_0^T\!\!\! \int_\S | u_\S(t,x)  - u_D(t,x)|^2 \, dt dx, 
$$
where $u_\S(t,x)$  is depending on $\O$ through the definition of the sets $\mathcal{O}$ and $\mathcal{S} = \Omega \setminus \overline{\mathcal{O}}$.

\smallskip
Notice, that the proof of Theorem~\ref{ThmMain} can be adapted to calculate the gradient of this type of functionals. The derivative is then calculated by replacing $U_M$ by $u_D$ in all the formulas in Theorem~\ref{ThmMain}. Consequently, the right hand side of the adjoint state problem \eqref{Eq-O-Adj} is modified and replaced by $u_\S - u_D$ and the first part of the shape derivative in \eqref{Grad} is replaced by $( u_\S - u_D)^2$. 

\smallskip
The gradient is initialized with an initial geometry with the disc at the bottom of the square with center coordinates $x_0 = 0.5$ and $y_0 = 0.2$. At each step, first the state $u$ is calculated with the current geometrical layout by solving \eqref{Eq-O}, next the adjoint state $g$ is obtained by solving \eqref{Eq-O-Adj} using the state data in the right hand side. Both problems are solved using the same values for $\kappa$, $R$, $U_M$ and $T$ used to calculate the temperature $u_D$ with the desired geometrical layout. Then the shape derivative is calculated using  the formula \eqref{Grad-ball} and $\Delta u_\S = \partial_t u_\S$ and $\Delta u_\O = \partial_t u_\O$ on the interface. This choice reduces the error on the shape derivative by replacing second order space derivatives with first order time derivative, which is approximated with a first order Euler scheme. The shape derivative is then used to update the current values of the coordinates of the center of the disc with respect to the constraints that require the disc to stay within the square domain. 

\begin{figure}[h!]
  \centering
  \makebox[\textwidth][c]{
  \begin{tabular}{lll}
    \includegraphics[width=0.3\linewidth]{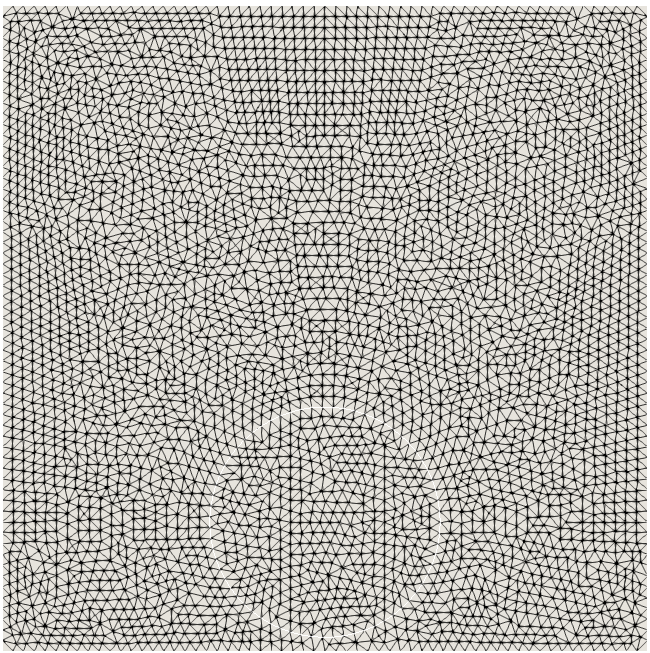}&
    \includegraphics[width=0.3\linewidth]{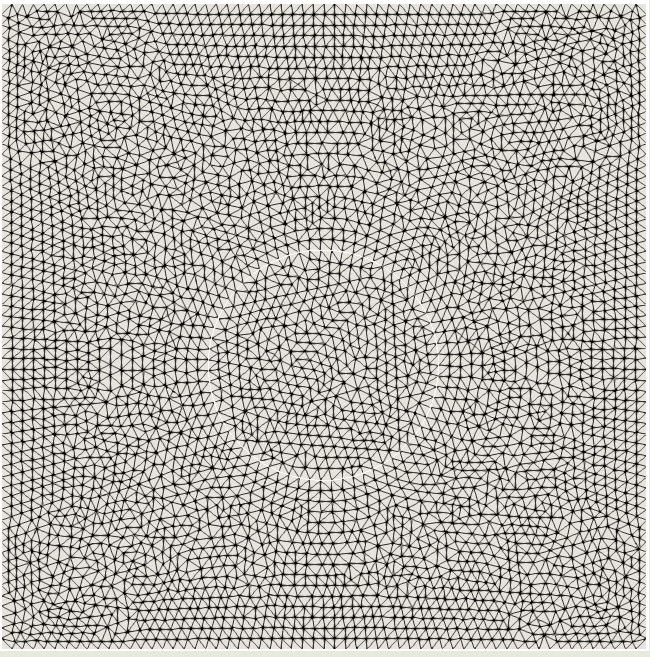}&
    \includegraphics[width=0.3\linewidth]{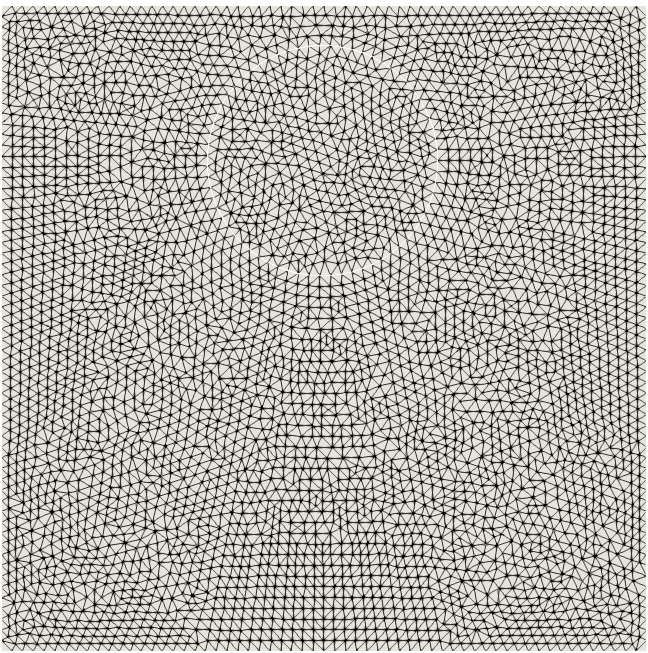}\\
  \end{tabular}
 }
\caption{Initial (left), intermediate (center) and final geometry (right) that corresponds to the desired geometry. \label{Fig:EvolCible}}
\end{figure}

\smallskip
The progress of the positions of the disc can be seen in Figure~\ref{Fig:EvolCible} with the initial position of the disc in the left figure, an intermediate position in the middle figure and the position at convergence in the right figure which corresponds to the desired position of the disc. 

\begin{figure}[h!]
  \centering
  \makebox[\textwidth][c]{
  \begin{tabular}{ll}
    \includegraphics[width=0.5\linewidth]{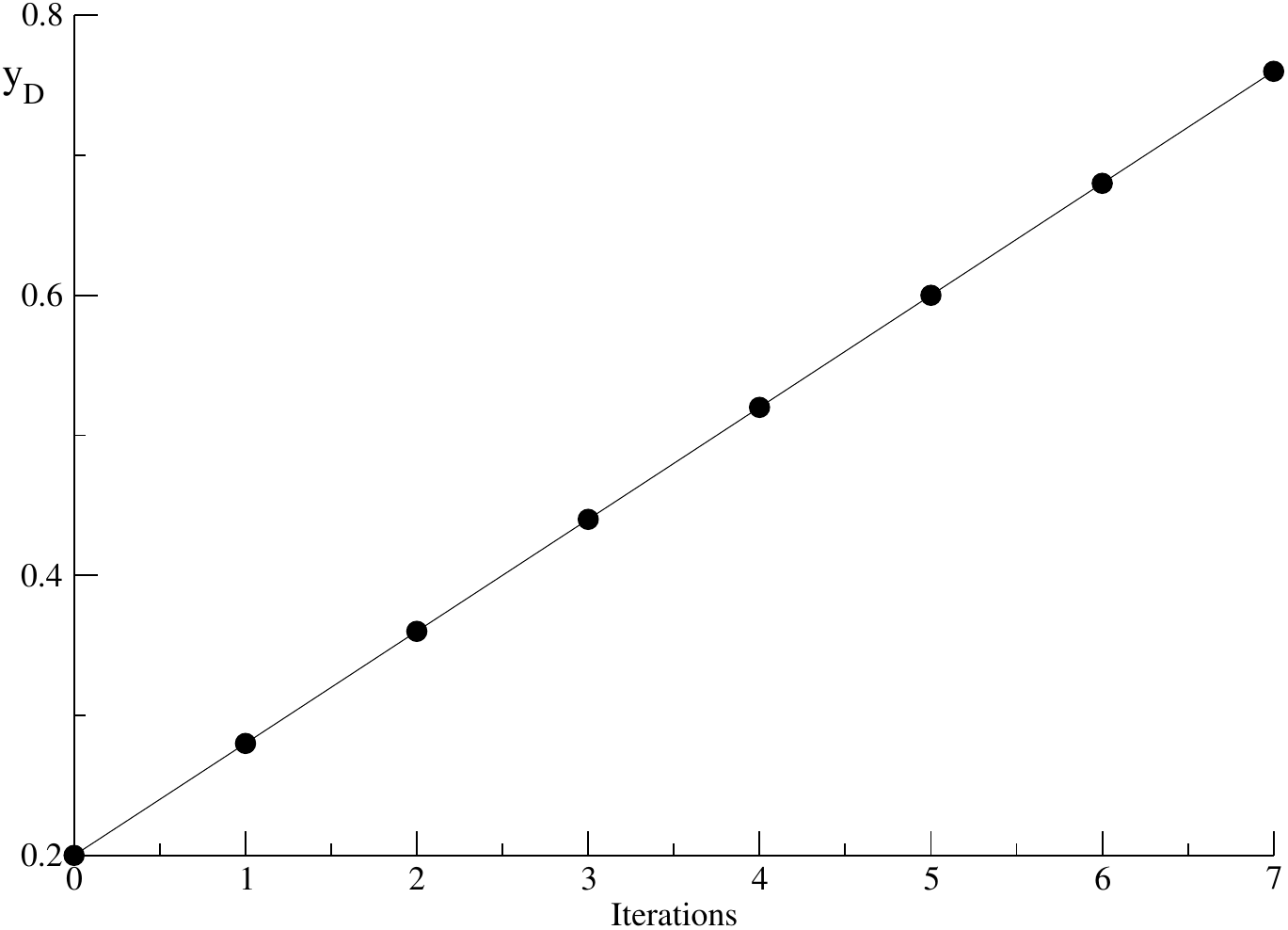}&
    \includegraphics[width=0.5\linewidth]{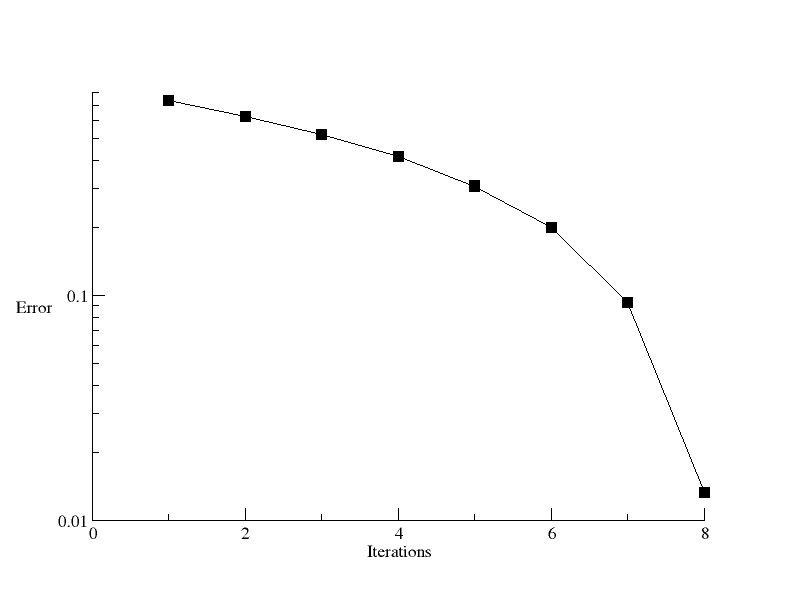}\\
  \end{tabular}
 }
\caption{Values of $y$ component of the disc's center versus number of iterations (left), error on the desired position at each iteration of the gradient (right). }\label{Fig:Ydisire}
\end{figure}

The choice of the initial position with $x_0 = x_D$ requires the gradient to move the position in the $y$ direction only. This is reported in Figure~\ref{Fig:Ydisire} where in the left side are presented the values of the $y$ component of the disc's center with respect to the number of iterations and in the right side is the error on the desired position at each iteration of the gradient. We can see that the gradient based on the shape derivative presented in this paper is very effective and allows to obtain a desired geometry knowing the associated temperature field solution to \eqref{Eq-O} with very few iterations (less than 10) and with fairly good accuracy with an error at convergence around $10^{-2}$.\\ 

\bigskip
\subsection{An optimization test case}
We consider the heat conduction problem described by equation \eqref{Eq-O} where $\O$ is a disc and $U_M$ is a known 
constant temperature imposed at the lower side of the squared domain $\Omega$ containing $\O$. The values for the thermal conductivity $\kappa$, the thermal contact resistance $R$ and the simulation time $T$ are taken the same as in the above test case. We consider solving the optimization problem consisting of identifying the optimal position of the disc $\O$ inside the square that allows to minimize the fitness function $J(\O)$ given in \eqref{Jfitness},  i.e.

\begin{equation*}
	J(\O) = \int_0^T \!\!\!\int_\S | u_\S(t,x)  - U_M|^2 \, dt dx, 
\end{equation*}
where $u_\S$ satisfies \eqref{Eq-O} is depending on $\O$ through the definition of the sets $\mathcal{O}$ and $\mathcal{S} = \Omega \setminus \overline{\mathcal{O}}$. Note that the Dirichlet boundary condition on the temperature in \eqref{Eq-O} is identical to the positive constant temperature $U_M$ taken in the definition of $J(\O)$. 

\medskip
In contrast to the validation test case, the optimal position of the disc is not known in advance. The optimization process is the same as in the validation test case described in the previous section, where at each iteration we first solve the equation \eqref{Eq-O}, and then the adjoint state equation~\eqref{Eq-O-Adj}. Next, the shape derivative is calculated using the formula \eqref{Grad-ball} from Section \ref{Subsec-Computation-Ball} and is used to update the coordinates of the disc's center. As in the validation case, the approximation $\Delta u_\S \approx \partial_t u_\S$ and $\Delta u_\O \approx \partial_t u_\O$ is used. The fitness function is reevaluated using \eqref{Jfitness} and the process is repeated several times. 
It is also noted that  we have taken the length of the square is greater than the diameter of the disc which is kept constant and corresponds to a volume of the disc equals to $10\%$ of the total volume. 

\begin{figure}[h!]
  \centering{
  \begin{tabular}{cc}
    \includegraphics[width=0.3\linewidth]{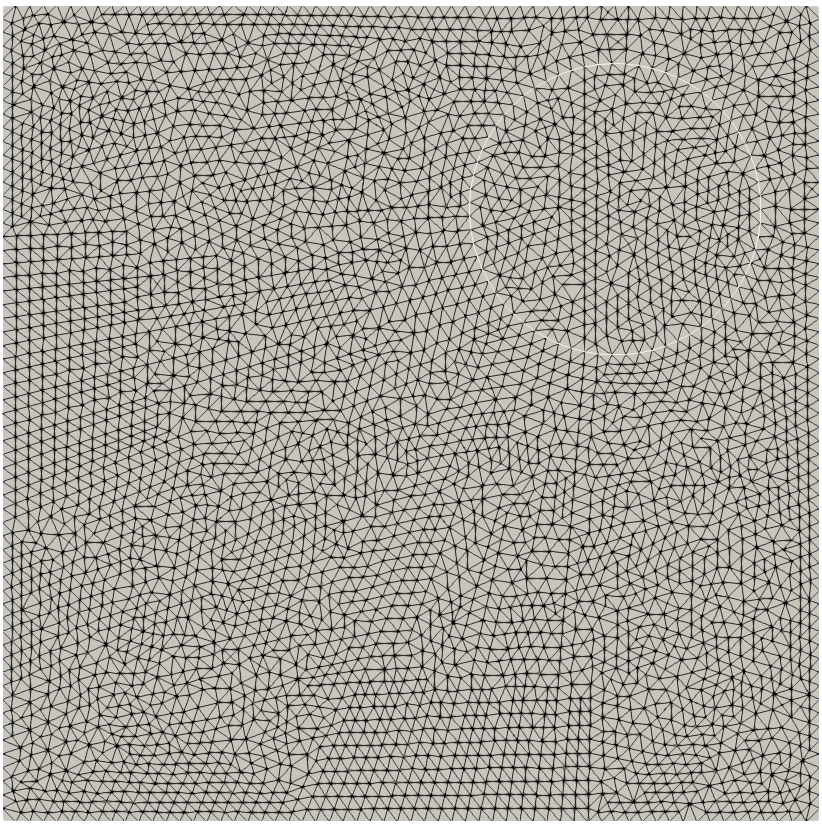}&
    \includegraphics[width=0.3\linewidth]{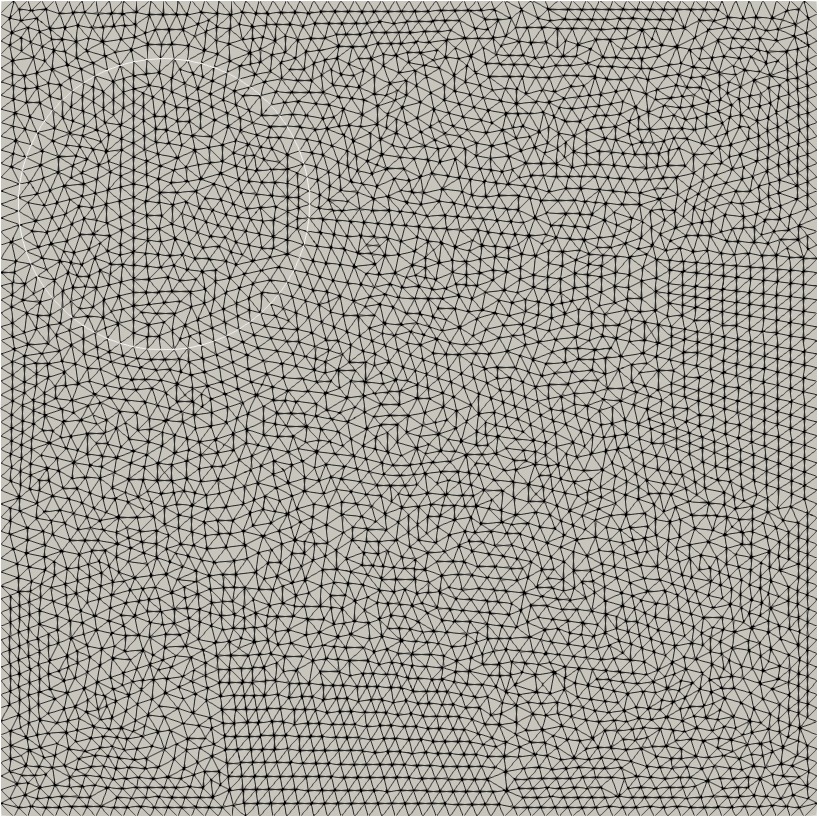}\\
    \includegraphics[width=0.3\linewidth]{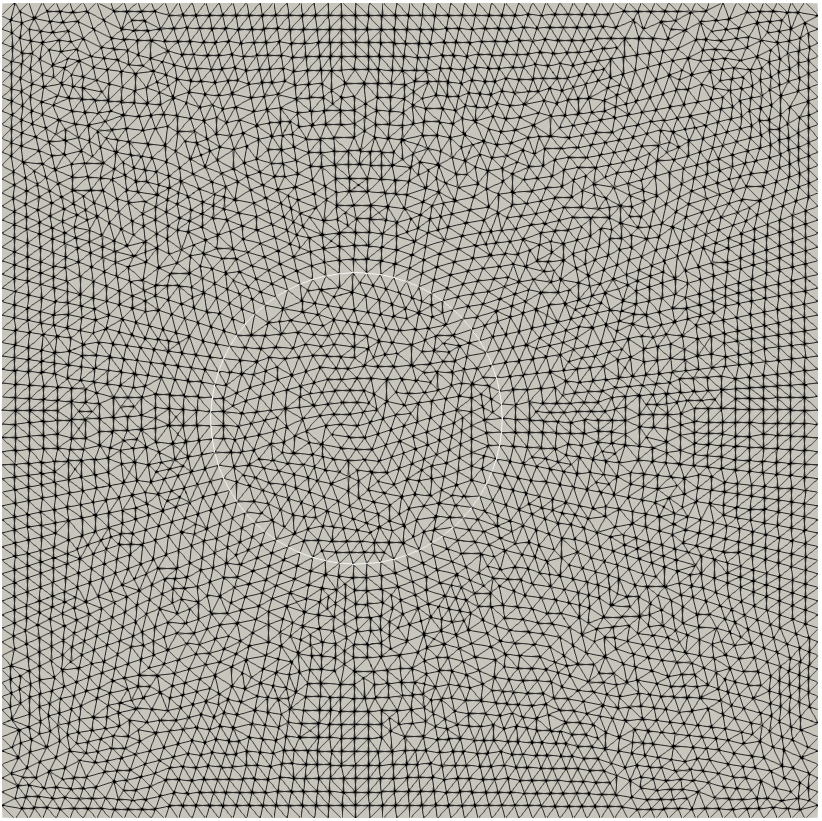}&
    \includegraphics[width=0.3\linewidth]{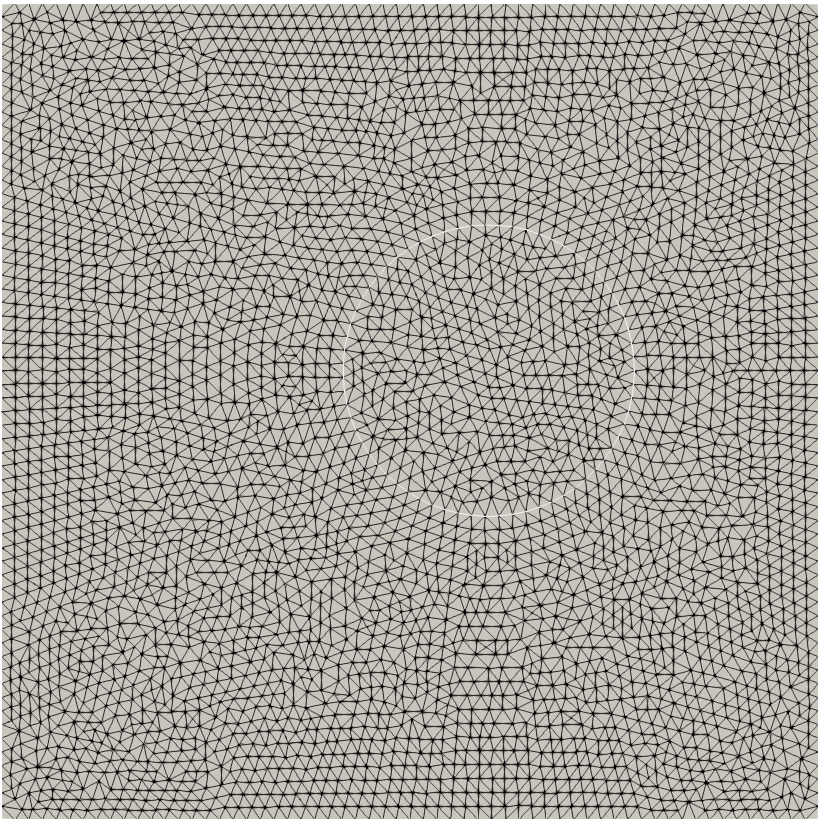}\\
    \includegraphics[width=0.3\linewidth]{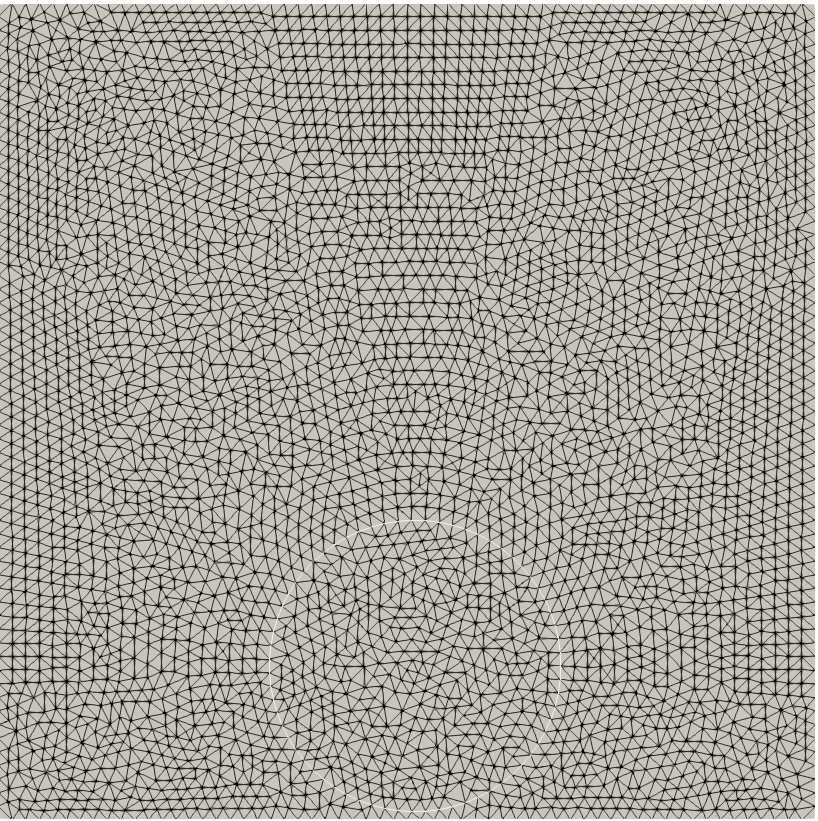}&
    \includegraphics[width=0.3\linewidth]{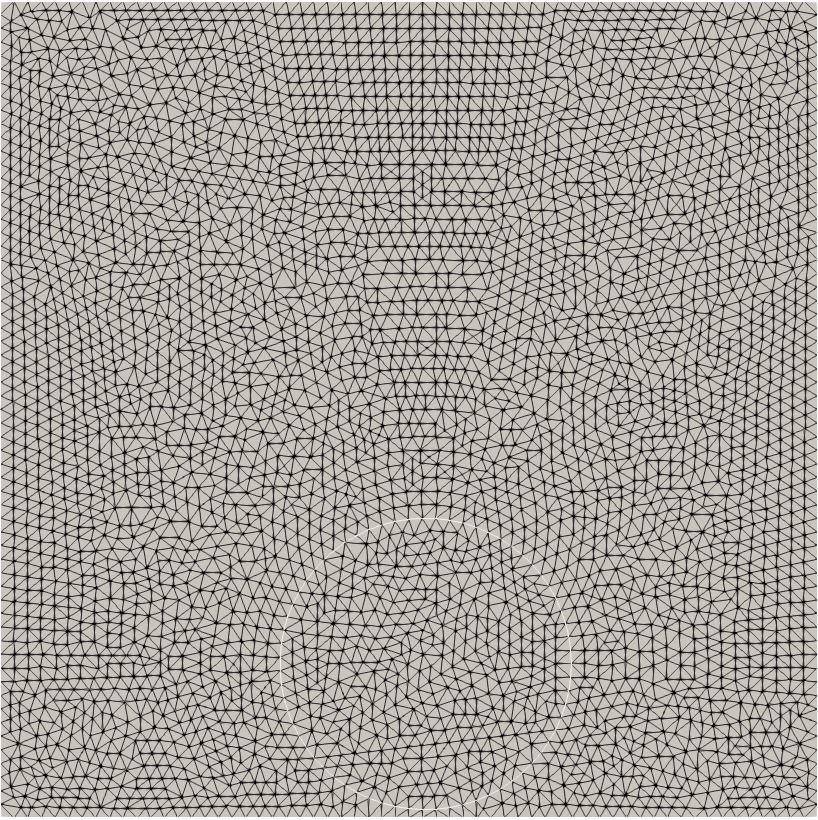}\\
  \end{tabular}
 }
\caption{From top to bottom, initial, intermediate and final solutions. Left, initial disc in the top right corner. Right, initial disc in the top left corner. \label{Fig:EvolUm1}}
\end{figure}

\medskip
We plot in Figure \ref{Fig:EvolUm1} the iterations corresponding to two symmetric initial positions for the disc, one on the top left corner, the other on the top right corner. Both test converge to the disc localized on the bottom and centered. This symmetry could have been expected in view of the symmetry of the problem, but we do not see any theoretical reasons to justify this property rigorously. Let us also notice that there is no guarantee that our method converges to a global optimal/minimizer, since we have not investigated the existence of a global optimal. However, in the ball class, we are minimizing a continuous function dependent on parameters (position) in a compact set, and so there is a global minimum, even if we do not know a priori whether our method converges to this global minimum. Nevertheless, when several distinct configurations all converge to the same thing, it is reasonable to believe that it is a global minimum. 

The obtained numerical results attest of the capabilities of the proposed method that allowed to obtain an optimal solution regardless of the initial guess. They also suggest that an optimal geometrical situation will correspond to the disc $\O$ closer to the border of the square at which the Dirichlet boundary condition is imposed. 

\medskip

\subsection{Mean temperature minimization case }
For this third test we consider the same heat conduction problem as in the first test namely the heat conduction problem described by equation \eqref{Eq-O} where $\O$ is a disc and $U_M$ is a known constant temperature imposed at the lower side of the squared domain $\Omega$ containing $\O$. The values for the thermal conductivity $\kappa$, the thermal contact resistance $R$ and the simulation time $T$ are taken the same as in the above test cases. We consider solving the minimization problem consisting of obtaining the optimal position of the disc $\O$ inside the square to minimize the quantity 
\begin{equation}
	\label{Def-J-3rd-test}
    J(\O) = \int_0^T \int_\S | u_\S(t,x)|^2 \, dt dx, 
\end{equation}
where $u_\S$ satisfies \eqref{Eq-O} is depending on $\O$ through the definition of the sets $\mathcal{O}$. This can be seen as the mean temperature in $\S$ over time. 

Again, the computation of the derivative of the function $J$ in \eqref{Def-J-3rd-test} can be done as in the proof of Theorem \ref{ThmMain}, as we have mentioned in Section~\ref{sec.numvc}. 

\begin{figure}[h!]
  \centering{
  \begin{tabular}{ccc}
    \includegraphics[width=0.33\linewidth]{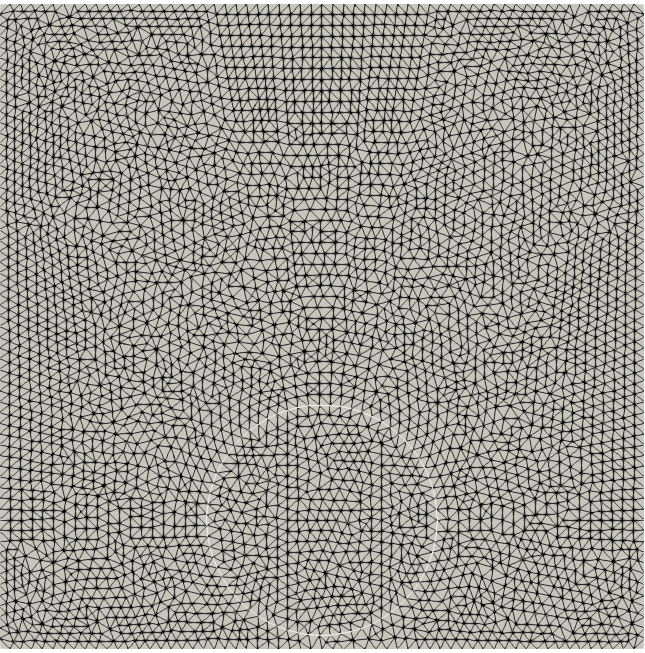}&
    \includegraphics[width=0.33\linewidth]{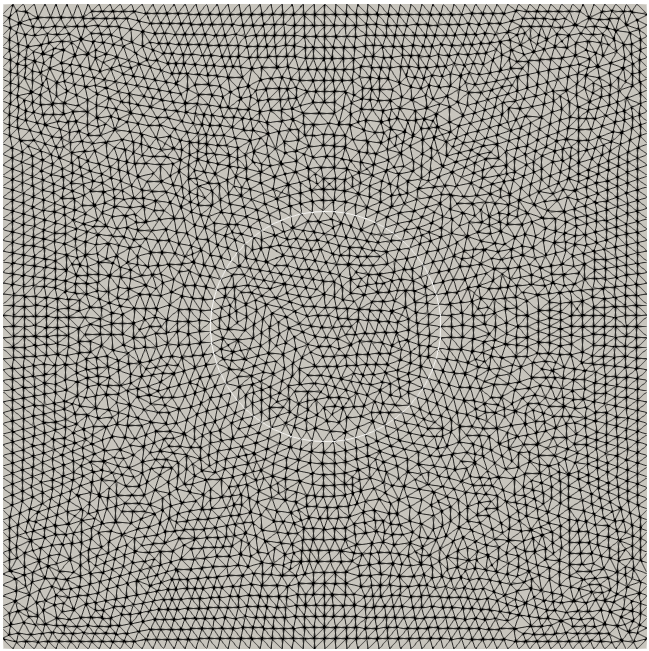}&
    \includegraphics[width=0.33\linewidth]{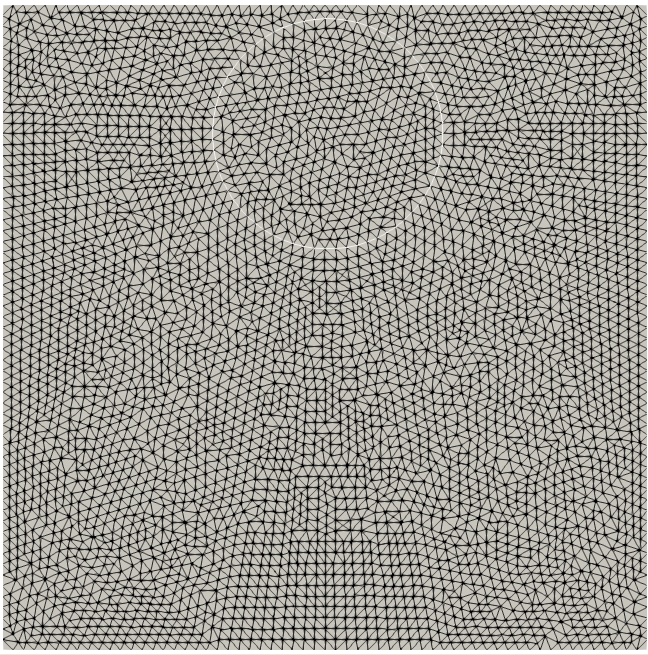}\\
  \end{tabular}
 }
\caption{Minimization of the mean temperature in $\S$ overtime. Initial, intermediate and final solutions. \label{Fig:EvolUm0}}
\end{figure}

The optimization process is the same as the above two tests with the initial position of the disc is at the bottom of the square. The results are given in Figure~\ref{Fig:EvolUm0}. We observe that the optimal geometrical layout is the one in which the disc is the farther from the heat source. This is quite reasonable since the disc contains a very conductive material and the closer it is from the heat source, the faster it will be heated and in turns increase the temperature in $\S$, which is the opposite of what we want.

%
\end{document}